\documentclass[11pt,a4paper,reqno,twoside]{amsart}
\usepackage{hyperref}

\usepackage{amssymb}
\usepackage{paralist}
\input xy
\xyoption{all}
\newtheorem{thm}{Theorem}[section]

\newtheorem*{prop*}{Proposition}
\newtheorem*{thm*}{Theorem}

\newtheorem*{lem*}{Lemma}

\newtheorem*{cor*}{Corollary}

\theoremstyle{definition}

\theoremstyle{remark}

\def\corner{$\hfill\lrcorner$}
\newtheorem{rema}[thm]{Remark}

\setdefaultenum{i.}{1.}{(a)}.
\newlength{\equwidth}
\def\equghost{\hspace{\equwidth}\ }

\def\A{\mathcal A}
\def\X{\mathfrak X}
\def\al{\alpha}

\def\ga{\gamma}
\def\de{\delta}

\def\ka{\kappa}
\def\kaf{{\kappa_f}}
\def\la{\lambda}

\def\om{\omega}
\def\Ga{\Gamma}
\def\De{\Delta}

\def\La{\Lambda}

\def\Om{\Omega}

\def\t{\otimes}

\def\goesto{\rightarrow}

\def\embed{\hookrightarrow}

\def\sgn{\mathrm{sgn}}

\def\G{\mathcal{G}}

\def\L{\mathcal{L}}

\def\ad{\mathrm{ad}}

\def\Ad{\mathrm{Ad}}

\def\Im{\mathrm{Im}}

\def\ddtz{{\frac{d}{dt}}_{|t=0}}

\def\Rh{{\/\hat{R}}}

\def\rr{\ensuremath{\mathbb{R}}}

\def\one{\ensuremath{\mathbb{I}}}
\def\wz{\ensuremath{\backslash\{0\}}}

\def\PO{\ensuremath{\mathrm{PO}}}

\def\O{\ensuremath{\mathrm{O}}}

\def\Aut{\ensuremath{\mathrm{Aut}}}

\def\inf{\ensuremath{\mathfrak{inf}}}

\def\co{\ensuremath{\mathfrak{co}}}

\def\o{\ensuremath{\mathfrak{so}}}
\def\so{\ensuremath{\mathfrak{so}}}

\def\h{\ensuremath{\mathfrak{h}}}
\def\k{\ensuremath{\mathfrak{k}}}

\def\g{\ensuremath{\mathfrak{g}}}
\def\n{\ensuremath{\mathfrak{n}}}

\def\p{\ensuremath{\mathfrak{p}}}
\def\gl{\ensuremath{\mathfrak{gl}}}

\def\Hol{\mathrm{Hol}}
\def\hol{\mathfrak{hol}}

\title{Homogeneous Cartan Geometries}
\author{Matthias Hammerl}
\address{Institut f\"ur Mathematik, Universit\"at Wien, Nordbergstra\ss e~15, A--1090 Wien, Austria}
\date{July 2007 (Some formulas corrected: May 2011)}

\email{matthias.hammerl@univie.ac.at}


\thanks{The author's work was supported by the IK I008-N funded by the University of Vienna}

\begin{document}

\begin{abstract}
  We describe invariant principal and Cartan connections
on homogeneous principal bundles and show how to calculate the curvature
and the holonomy; in the case of an invariant Cartan connection
we give a formula for the infinitesimal automorphisms.
The main result of this paper is that
the above calculations are purely algorithmic. As an example
of an homogeneous parabolic geometry we treat a conformal structure
on the product of two spheres.
\end{abstract}

\maketitle

\subsection*{Introduction}
We begin with a discussion of invariant principal connections
in section \ref{section-invpri}\ and present a similar
treatment of invariant Cartan connections in section
\ref{section-homcart}. For both kinds of connections we
give an explicit description of the holonomy Lie algebra
which is due to H.C. Wang.
Using ideas from \cite{cap-infinitaut}
 we can use the knowledge of the holonomy Lie algebra to calculate
the infinitesimal automorphisms of a Cartan geometry in
section \ref{subsection-infaut}.

Nothing in here is really new. Invariant
principal connections were already treated by
 H.C. Wang in \cite{wang58}. The case of invariant
Cartan connections is quite analogous.
Infinitesimal automorphisms of Cartan geometries
were discussed by A. \v{C}ap in \cite{cap-infinitaut} and the consequences drawn here for the homogeneous case
are quite elementary.

Nevertheless, as a whole, this provides a nice and simple framework for
treating homogeneous parabolic geometries: for
such structures, among which are conformal
structures, almost CR structures of hypersurface type, projective structures, projective contact structures and almost quaternionic
structures, one has an equivalence of categories with
certain types of Cartan geometries (\cite{tanaka79,yamaguchi,cap-schichl,cap-slovak-par}): this immediately extends
all notions discussed for Cartan geometries, like curvature, holonomy and infinitesimal
automorphisms, to these structures.

In section \ref{section-holsph}\ we consider
as an explicit example the case of a conformal structure
on the product of two spheres: We use a well known method to obtain
the canonical Cartan connection to this parabolic geometry and then
calculate the holonomy Lie algebra of this geometry.

\subsubsection*{\textbf{Acknowledgments}}
Special thanks go to Andreas \v{C}ap for his support and encouragement.
I also thank Katja Sagerschnig for her valuable comments.

\newpage

\def\ka{\kappa}
\def\kaf{\kappa}

   \section{Invariant principal connections}\label{section-invpri}

   \subsection{Homogeneous $P$-principal bundles}\label{section-homogeneous-bundles}
   Let $\G\overset{\pi}{\goesto} M$\ be a $P$-principal bundle.
   The right
   action of $P$\ on $\G$\ shall be called
   \begin{align*}
     r:\G\times P&\goesto \G,\\
     r^p(u)&=r(u,p)=u\cdot p.
   \end{align*}
   We say that the $P$-principal bundle $\G$\ is homogeneous
   when there is a Lie group $H$\ acting fiber-transitively
   on $\G$\ by principal bundle automorphisms.

   We denote the action of an element $h\in H$\
   on $\G$\ by $\la_h\in\Aut(\G)$\, and simply write $\la_h(u)=h\cdot u$\
   for $u\in \G$. The diffeomorphism $\la_h$\ factorizes to a
   diffeomorphism $\check\la_h$\ of $M$\ and fiber-transitivity of
   the action of $H$\ on $\G$\ is equivalent to transitivity of the
   induced action of $H$\ on $M$.
   
   Let $o\in M$\ be an arbitrary point. Then
   the isotropy group $K=H_o$\ of $o$\ is a Lie subgroup of $H$\ and
   $M=H/K$. Now the action of $K$\ on $\G$\ leaves the fiber $\G_o$\
   of $\G$\ over $o$\ invariant. Take some $u_0\in\G_o$;
   then one sees that there is a unique homomorphism of
   Lie groups $\Psi:K\goesto P$\ with $k\cdot u_0=u_0\cdot \Psi(k)$.
   Now one checks that in fact $\G=H\times_K P=H\times_{\Psi} P$,
   the associated bundle to the principal $K$-bundle $H\goesto H/K$\
   obtained by the action of $K$\ on $P$\ by $\Psi$. We denote
   the equivalence class
   \begin{align*}
     {\tilde{\pi}}(h,p)=\{(h\cdot k,\Psi(k)\cdot p),k\in K\}=[h,p].
   \end{align*}

   So we described an arbitrary
homogeneous $P$-principal bundle as a quotient of the
trivial bundle $H\times P$. We have a $K$-principal bundle whose
base is a $P$-principal bundle:
\begin{align*}
  \xymatrix{
H\times P \ar[d]^-{{\tilde{\pi}}} & \ar[l] K \\
H\times_K P \ar[d]^-{\pi} & \ar[l] P\\
H/K
}
\end{align*}

\vspace{0.1cm}
   We summarize: Every homogeneous $P$-principal bundle is of the
   form $\G=H\times_K P\goesto H/K$ with the canonical left- respectively right-
   actions of $H$\ resp. $P$\ on $\G$. We can also say: the data defining
   an $H$-homogeneous $P$-principal bundle is a 4-tuple $(H,K,P,\Psi)$\
   , where $H$\ and $P$\ are Lie groups, $K$\ is a closed subgroup of $H$\
   and $\Psi$\ is a homomorphism of Lie groups from $K$\ to $P$.
   
   \subsection{Invariant principal connections on homogeneous $P$-principal bundles}
\label{subsection-invprinc}
   Any $H$-invariant $P$-principal connection $\ga\in\Om(H\times_K P,\p)$\
   can be pulled back to a $P$-principal connection on the trivial
   bundle $H\times P$.
   Now Frobenius reciprocity (see e.g. \cite{michor-topics}, 22.14) provides a one-to-one correspondence
   between $K$-invariant horizontal forms on $H\times P$\ and
   $\p$-valued forms on $H\times_K P$. 
   This can be used to show that $H$-left invariant $P$-principal connections on $H\times_K P$\ correspond exactly to
   certain linear maps $\al:\h\goesto\p$. Precisely:
   Left-trivialize $H\times P=H\times P\times\h\times\p$, then
   \begin{thm}[\cite{wang58}, Prop. (5.1), Prop. (5.2);
       \cite{mrh-dip}, Thm. 2.2.6, Thm.4.1.1;]
     Every invariant principal connection $\ga$\ on $\G=H\times_K P$\
     is obtained by factorizing a $\p$-valued one form
     \begin{align*}
       \hat\ga\in\Om^1_K(H\times P,\p)^{hor},\\
       (h,p,X,Y)\mapsto \Ad(p^{-1})\al(X)+Y
     \end{align*}
     with $h\in H, p\in P,X\in \h$\ and $Y\in\p$.
     Here the conditions on $\al\in L(\h,\p)$\ such that $\hat\ga$\
     is indeed $K$-invariant and horizontal are
     \begin{enumerate}
     \item $\al_{|\k}=\Psi'$
     \item $\al(\Ad(k)X)=\Ad(\Psi(k))\al(X)$.
     \end{enumerate}
     I.e.: $\al$\ is a $K$-equivariant extension of $\Psi'$\ to
     a linear map from $\h\goesto \p$.
   \end{thm}

   This description of invariant principal connections on
   $\G$\ leads to an easy formula for the curvature. If
   $\ga\in\Om^1(\G,\p)$\ is the invariant principal connection
   corresponding to $\al:\h\goesto\p$,
   the curvature ${\bar\rho}\in\Om^2(\G,\p)$\ of $\ga$\ is given by
   \begin{align*}
{\bar\rho}(\xi,\eta)=-\ga([\xi_{hor},\eta_{hor}]),
   \end{align*}
   where $\xi_{hor},\eta_{hor}$\ are the horizontal projections
   of vector fields $\xi,\eta\in\X(\G)$.

   Now ${\bar\rho}$\ is $P$-equivariant and horizontal and thus
   factorizes to a $H\times_K \p$-valued 2-form  on $H/K$,
   i.e., a section of $H\times_K \La^2(\h/\k)\t\p$.
   It is easy to see that $H$-invariance of $\ga$\ implies $H$-invariance of
   this section, which is therefore determined by a unique
   $K$-invariant element $\rho$\ of $\La^2(\h/\k)\t\p$.
   
   Therefore we will say that the curvature of $\al$\
   is $\rho=\rho_{\al}\in\La^2(\h/\k)\t\p$,
   and one calculates that
   \begin{align*}
        \rho(X_1,X_2)=[\al(X_1),\al(X_2)]-\al([X_1,X_2])
   \end{align*}
   for $X_1,X_2\in\h$.
   
   So the curvature of $\al$\ is its failure to be (an extension
   of $\Psi':\k\goesto\p$\ to) a  homomorphism of Lie algebras
   $\h\goesto\p$.

   Denote the holonomy of the invariant principal connection corresponding
   to $\al:\h\goesto\p$\ by $\hol(\al)$.
   H.C. Wang gave an explicit description of $\hol(\al)$:
   \begin{thm}[\cite{wang58}, Theorem (B)]\label{homogeneous-holonomy}
     Denote by $\Rh=\langle\{\rho(X_1,X_2)|X_1,X_2\in\h\}\rangle$\
     the span of the image of $\rho$\ in $\p$.
     Then $\hol(\al)$\ is the $\h$-module generated by $\Rh$. I.e.:
     \begin{align}\label{holonomy-formula}
       \hol(\al)=\Rh+[\al(\h),\Rh]+[\al(\h),[\al(\h),\Rh]]+\cdots.
     \end{align}
   \end{thm}
   Note that it follows in particular 
   that $\Rh+[\al(\h),\Rh]+[\al(\h),[\al(\h),\Rh]]+\cdots$\
   is already a Lie subalgebra of $\p$.
   The proof is quite involved. From the Ambrose-Singer theorem
   (\cite{ambrose-singer-53}) one knows hat
   $\Rh\subset\hol(\al)$. The essential part is then the 
   construction of a group which can be explicitly described and which
   contains the holonomy group as a normal subgroup.
   Then $\hol(\al)$\ is a module under this group and is shown to be
   generated by $\Rh$, which can then be reformulated as \eqref{holonomy-formula}.
   
   \section{Homogeneous Cartan Geometries}\label{section-homcart}
   One has a very similar description of invariant Cartan connections
   on homogeneous principal bundles. Recall
   that for Lie groups $(G,P)$ with $P$\ a closed subgroup of $G$\
   a Cartan geometry of type $(G,P)$\ is a $P$-principal bundle $\G$\
   over a manifold $M$\ endowed with a one form
   $\om\in\Om^1(\G,\g)$, called the Cartan connection, satisfying
   \begin{enumerate}
   \item $\om$\ is $P$-equivariant:\ $(r^p)^*(\om)=\Ad(p^{-1})\circ\om$\ for all $p\in P$.
   \item $\om$\ reproduces fundamental vector fields:\
     $\om(\ddtz u\cdot \exp(tY))=Y$\ for all $Y\in\p$.
   \item $\om$\ is an absolute parallelism:\ at every $u\in \G$\
     the map $\om_u:T_u\G\goesto \g$\ is an isomorphism.
   \end{enumerate}

   \noindent
   We will say that $(\G,\om)$\ is a Cartan geometry of type $(G,P)$.

   A Cartan geometry $(\G,\om)$\ is homogeneous if there is a Lie group
   $H$\ acting fiber-transitively on $\G$\ by automorphisms of the
   Cartan geometry. I.e., when we denote the action of $H$\ on $\G$\
   by $\la_h:\G\goesto\G$\ we have that for every $h\in H$\
   \begin{enumerate}
   \item $\la_h$\ is an automorphism of the $P$-principal bundle $\G$:\ 
     $\la_h(u\cdot p)=\la_h(u)\cdot p$
   \item $\la_h$\ preserves $\om$:\ $\la_h^*(\om)=\om$.
   \end{enumerate}

   Now we saw in section \ref{section-homogeneous-bundles}\
   that $\G$\ is of the form $H\times_{\Psi} P$\ for
   some homomorphism $\Psi:K\goesto P$\ and analogously
   to the case of principal connections one gets the
   following description of invariant Cartan connections of
   type $(G,P)$\ on $\G=H\times_K P$:
   \begin{thm}[\cite{wang58}, Thm. 4; \cite{mrh-dip}, Thm. 4.2.1.]
     Invariant Cartan connections on $H\times_K P$\ are in
     1:1-correspondence with maps $\al:\h\goesto\g$\
     satisfying
     \begin{enumerate}[(C.1)]
     \item \label{cartcon1} $\al_{|\k}=\Psi'$ 
     \item \label{cartcon2}$\al(\Ad(k)X)=\Ad(\Psi(k))\al(X)$\ for all $X\in\h,k\in K$
     \item \label{cartcon3}$\al$\ induces an isomorphism of $\h/\k$\ with $\g/\p$.
     \end{enumerate}
     Explicitly: given such an $\al$, the corresponding Cartan connection
     $\om$\ is obtained by factorizing
     \begin{align*}
       &\hat\om\in\Om^1(H\times_K P,\g),\\
       &\hat\om((h,p,X,Y))=\Ad(p^{-1})\al(X)+Y.
     \end{align*}
   \end{thm}

   \noindent
   We will say that $\al$\ is a Cartan connection.

   Similarly as in the case of invariant principal connections,
   the curvature of an invariant Cartan connection $\al:\h\goesto\p$\
   is described by an element $\ka=\ka_{\al}\in\La^2(\h/\k)\t\g$:
   it is given by
   \begin{align*}
     \ka(X_1,X_2)=[\al(X_1),\al(X_2)]-\al([X_1,X_2].
   \end{align*}
   Using the isomorphism induced by $\al$\ between
   $\h/\k$\ and $\g/\p$ we can also regard the curvature as
   \begin{align*}
     \ka=\ka(\al)\in \La^2(\g/\p)^*\t\g.
   \end{align*}

   \subsection{The holonomy of a Cartan connection}\label{holofcart}
   Consider a Cartan geometry $(\G,\om)$\ of type $(G,P)$.
   We can extend the structure group of $\G$\ from
   $P$\ to $G$\ by taking $\G'=\G\times_P G$.
   Then a Cartan connection $\om\in\Om^1(\G,\g)$\ extends
   equivariantly to a $G$-principal connection $\om'\in\Om^1(\G',\g)$\
   on $\G'$.
   We will say that the holonomy Lie group $\Hol(\om)$\ of the Cartan geometry
   $(\G,\om)$\ is $\Hol(\om')$.

   When $(\G,\om)$\ is homogeneous, i.e., $\G=H\times_K P$\
   and $\om$\ is induced by an $\al:\h\goesto\g$\ satisfying
   (C.\ref{cartcon1})-(C.\ref{cartcon3}) the above construction yields
   $\G'=H\times_{\Psi} G$, with $\Psi$\ regarded as a map from
   $K\goesto P\embed G$, and $\om'$\ corresponds (again) to
   $\al:\h\goesto\g$; Note here that $\al$ in particular satisfies (C.\ref{cartcon1})\
   and (C.\ref{cartcon2}).
   
   Thus we have
   \begin{thm}\label{hol-theorem}
     The holonomy Lie algebra of a homogeneous Cartan geometry \\
     $(H\times_K P,\al:\h\goesto\g)$\ is
     \begin{align*}
       \hol(\al)=\Rh+[\al(\h),\Rh]+[\al(\h),[\al(\h),\Rh]]+\cdots,
     \end{align*}
     where $\Rh=\langle\{\ka(X_1,X_2)|X_1,X_2\in\h\}\rangle$.
   \end{thm}

\subsection{Infinitesimal Automorphisms}\label{subsection-infaut}
Let $(\G\goesto M,\om)$\ be a Cartan geometry of
type $(G,P)$. A vector field $\xi\in\X(\G)$\ is
an \emph{infinitesimal automorphism}\ of $(\G\goesto M,\om)$\
if it satisfies
\begin{enumerate}[(1)]
\item \label{autcond1} $\xi$\ is $P$-invariant:\ $Tr^p \xi(u)=\xi(u\cdot p)\ 
  \forall u\in\G,p\in P$
\item \label{autcond2} $\xi$\ preserves $\om$:\ $\L_{\xi}\om=0$.
\end{enumerate}

We remark that the Lie algebra of the automorphism group of
the Cartan geometry $(\G,\om)$ is formed by the  \emph{complete}\ vector
fields on $\G$\ which are infinitesimal automorphisms.

Condition (\ref{autcond1})\ is equivalent to
$\om\circ \xi:\G\goesto\g$\ being $P$-equivariant.
Thus an infinitesimal automorphism $\xi$\ of $(\G,\om)$\
gives rise to a section of the adjoint tractor bundle
\begin{align*}
  \A M:=\G\times_P \g.
\end{align*}

We want to describe, in terms of geometric data on $\A M$, those sections of $\A M$\ which correspond to infinitesimal automorphisms of $(\G,\om)$. We first show how $\om$\ induces a linear connection on $\A M$:
Since the action of $P$\ on $\g$\ is just the restriction of the adjoint
action of $G$\ on $\g$\ we have
\begin{align*}
  \A M=\G\times_P\g=(\G\times_P G)\times_G\g=\G'\times_G\g.
\end{align*}
Recall from \ref{holofcart}: $\G'=\G\times_P G$: $\G'$\ is the
extension of structure group of $\G$\ from $P$\ to $G$\ and
$(\G',\om')$\ is a $G$-principal bundle endowed with a principal
connection.
Thus the principal connection $\om'$\ on $\G'$\ induces
a linear connection $\nabla$\ on $\A M$.

The second ingredient we need to give a condition on a section $s\in\Ga(\A M)$\
to be an automorphism of $(\G,\om)$\ comes from the curvature
function
\[\kaf:\G\goesto \La^2(\g/\p)^*\t \g.\]
It is a $P$-equivariant smooth map, and since $TM=\G\times_P\g/\p$,
\begin{align*}
  \kaf\in\Ga(\La^2(T^*M)\t\A M)\subset\Ga(T^*M\t T^*M\t \A M).
\end{align*}
Now note that since $\A M=\G\times_P \g$\ the canonical surjection
$\Pi:\g\goesto \g/\p$\ induces a map
\begin{align*} 
\bar\Pi:\A M\goesto TM.
\end{align*}
For $s\in\Ga(\A M))$\ we may thus consider
\begin{align*}
  i_s\kaf:=\kaf(\bar\Pi(s),\cdot)\in\Om^1(M,\A M).
\end{align*}

The following theorem characterizes infinitesimal automorphisms
of a parabolic geometry as parallel sections of a connection
on the adjoint tractor bundle $\A M$:
\begin{thm}[\cite{cap-infinitaut}, Prop. 3.2]\label{thm-inf}
  A section $s$\ of $\A M$\ corresponds to an infinitesimal automorphism
  of $(\G,\om)$\ if and only if 
\begin{align*}    
\nabla s+i_s\kaf=0.
\end{align*}
  \begin{proof}
    Let $\xi\in\X(\G)$\ be a $P$-invariant vector field on $\G$.
    Then $\om\circ\xi$\ is a $P$-equivariant map $\G\goesto\g$\
    and the corresponding section $s\in\Ga(\A M)$\ is obtained
    by factorizing $u\mapsto [u,\om(\xi(u))]$.
    Now $\L_\xi\om=i_{\xi}d \om+d(\om(\xi))$;
    Take a vector field $\eta\in\X(\G)$\ which is 
    $\pi$-related to a vector field $\tilde\eta$\ on $M$; Then
    \begin{align*}
      (\L_\xi\om)(\eta)=d\om(\xi,\eta)+\eta\cdot\om(\xi)
      =\ka(\xi,\eta)+\eta\cdot\om(\xi)-[\om(\xi),\om(\eta)].
    \end{align*}
    Since $\nabla_{\tilde\eta} s$\ corresponds to
    to the $P$-equivariant map
    $u\mapsto \eta(u)\cdot\om(\xi)+\ad_{\om_u(\eta)}(\om_u(\xi))$
    this proves the claim.
  \end{proof}
\end{thm}

Now
\begin{align}\label{definition-hatnabla}
  \hat\nabla_{\xi}s&:=\nabla_{\xi}s+\kaf(\bar\Pi(s),\xi)\
  \mathrm{for}\ \xi\in\Ga(TM),s\in\Ga(\A M)
\end{align}
is again a linear connection on $\A M$,
and thus Theorem \ref{thm-inf}\ says that
the infinitesimal automorphisms of $(\G,\om)$\
are
\begin{align*}
  \inf(\om):=\{s\in\Ga(\A M):\ \hat\nabla s=0\}.
\end{align*}
i.e.:$\inf(\om)$\ consists of the parallel sections of $(\A M,\hat\nabla)$.

This allows us to reformulate the problem of determining $\inf(\om)$\
in the following way:
\begin{thm}\label{thm-autom}
  \begin{align*}
    \inf(\om)=\{X\in\g:\Ad(\Hol(\hat\nabla))X=\{X\}\}.
  \end{align*}
  In particular, if $M$\ is simply connected,
  \begin{align*}
    \inf(\om)=\{X\in\g:\ad(\hol(\hat\nabla))X=\{0\}\}.
  \end{align*}
\end{thm}
This follows from the well known fact that parallel sections
of vector bundles correspond to holonomy-invariant elements of
the modeling vector space.

\subsubsection{Infinitesimal Automorphisms of Homogeneous Cartan Geometries}
We can now apply theorem \ref{thm-autom}\ to the case of a
Cartan connection $\al:\h\goesto\g$\ on a homogeneous principal
bundle $H\times_K P\goesto H/K$\ as discussed in section \ref{section-homcart}:

$\hat\nabla$\ as defined in \eqref{definition-hatnabla}\
is $H$-invariant and it is easy to see that it
is induced by the $G$-principal connection
\begin{align*}
  \hat\al&:\h\goesto\gl(\g),\\
  \hat\al(X)&=\ad(\al(X))+\kaf(\Pi(\cdot),\al(X)+\p)
\end{align*}
on $G'=\G\times_P G$.

Thus we have
\begin{thm}
  Let $\al:\h\goesto\g$\ be a Cartan connection of type $(G,P)$\
  on a simply connected homogeneous space $H/K$.
  Then the Lie algebra of infinitesimal automorphisms of the corresponding
  homogeneous Cartan geometry consists of all elements of
  $\g$\ which are stabilized by $\hol(\hat\al)$, i.e.:
  \begin{align*}
    \inf(\al)=\{X\in\g:\hol(\hat\al)X=\{0\}\}.
  \end{align*}
\end{thm}
Since we can determine $\hol(\hat\al)$\ by using theorem \ref{homogeneous-holonomy}\ it is a purely algorithmic task to calculate the infinitesimal
automorphisms of a homogeneous Cartan geometry.

Of course we know that $H$\ acts by automorphisms of Cartan geometries
on $(\G,\al)=(H\times_K P,\al)$\ from the left. It is clear that the fundamental vector
fields on $\G$\ for this action are infinitesimal automorphisms
and we have seen above that they are thus determined by elements of
$\g$:
It is easy to see
that these elements are exactly those in the image of
$\al:\h\goesto\g$.
It is not difficult either to verify that indeed
$\al(\h)$\ lies in the kernel of every element of
 the holonomy Lie algebra of $\hat \al$; the main observation
here is that $\hat\al(X)\al(Y)=\al([X,Y])$\ for $X,Y\in\h$.

   \section{The conformal holonomy of the product of two
     spheres}\label{section-holsph}
It is a classical result of \'{E}lie Cartan (\cite{cart-con}) that conformal geometries
are equivalent to certain parabolic geometries.
Thus \ref{holofcart}\ provides a notion of holonomy for conformal structures, which is
called \emph{conformal\ holonomy}. Conformal holonomies
induced by bi-invariant metrics on Lie groups have been treated in \cite{leitner-conformal}.
Our setting for calculating conformal holonomies works for invariant
conformal structures on arbitrary simply connected homogeneous spaces.

\def\gone{\ensuremath{{g_1}}}
\def\gtwo{\ensuremath{{g_2}}}
\def\hatgone{\ensuremath{{{\overline{g}}_1}}}
\def\hatgtwo{\ensuremath{{{\overline{g}}_2}}}
\def\hatgs{\ensuremath{\overline{g}_{(s,s')}}}
\def\deone{\ensuremath{{\de_1}}}
\def\detwo{\ensuremath{{\de_2}}}
\def\gs{\ensuremath{g_{(s,s')}}}
\def\tildka{\mbox{$\tilde\ka$}}

   Let $S^p,S^q$\ be the Euclidean spheres of dimension $p,q$\ with
   $p+q\geq 3$\ and let \hatgone,\hatgtwo\ denote their Riemannian
   metrics of radius $1$.  For $s\in\rr,s>0$\ and $s'\in\rr\wz$
   we have the (pseudo-)Riemannian metric
   $\hatgs=(\frac{1}{s}\hatgone,\frac{1}{s'} \hatgtwo)$\
   on $M=S^p\times S^q$. When $s'>0$\ $\hatgs$\ is
   positive definite and for $s'<0$\ it has signature $(p,q)$. The
   conformal class $[\hatgs]$\ of this Riemannian
   metric endows $S^p\times S^q$\
   with a conformal structure and we are going to calculate its
   conformal holonomy. To do this, we first switch to a homogeneous or Lie group
   description of the conformal geometry $(S^p\times S^q,\hatgs)$.

   Since $S^p=O(p+1)/O(p)$\ we have
   $M=S^p\times S^q=H/K$\ with $H=O(p+1)\times O(q+1)$\
   and $K=O(p)\times O(q)$.
   We will write elements of $\o(p+1)$\ as
   \begin{align*}
     v\oplus A=
     \left(
       \begin{matrix}
       0 & -v^t \\
       v & A         
       \end{matrix}
     \right)
   \end{align*}
   with $v\in\rr^p$\ and $A\in\o(p)$.
   So $\o(p+1)=\rr^p\oplus\o(p)$\ and the Lie bracket is
   \begin{align}\label{bracket-opplus1}
     [v_1\oplus A_1,v_2\oplus A_2]=\bigl(A_1 v_2-A_2 v_1\bigr)\oplus
     \bigl(v_2 v_1^t- v_1 v_2^t +[A_1,A_2]\bigr).
   \end{align}

   Thus $\k=\o(p)\oplus\o(q)$\ and $\h=(\rr^p\oplus\rr^q)\oplus\k$.
   We will denote $\n=\rr^p\oplus\rr^q<\h$.
   Let $\gone=\one_p$,\ $\gtwo=\one_q$\ be the standard Euclidean inner
   products on $\rr^p$\ and $\rr^q$. It will later be useful also
   to regard $\gone$,$\gtwo$\ as (degenerate) bilinear forms on $\rr^{p+q}$\ by
   trivial extension.

   It is easy to see that the
   $H=\O(p+1)\times \O(q+1)$-invariant (pseudo-)Riemannian metric
   $\hatgs=(\frac{1}{s}\hatgone,\frac{1}{s'}\hatgtwo)$\ on $M=S^p\times S^q=H/K$\
   corresponds to the $K$-invariant (pseudo-) inner product
   $\gs=\frac{1}{s}\gone\oplus \frac{1}{s'}\gtwo$\ on $\rr^p\oplus\rr^q$.

   \subsection{The prolongation to a canonical Cartan connection}

   Now we describe $(H/K=S^p\times S^q,[\hatgs])$\ as a
   (canonical) Cartan geometry of type $(G,P)$: here
   $G=\PO(p+q+1,1)$\ for $s'>0$\ and $G=\PO(p+1,q+1)$\ for
   $s'<0$.
   Let $g:=\gone+\sgn(s')\gtwo$. This is the standard inner product
   of the same signature as $\gs$.
   The Lie algebra of $G$\ is graded
   \begin{align*}
     \g=\g_{-1}\oplus\g_0\oplus \g_1=\rr^{p+q}\oplus\co(\rr^{p+q},g)\oplus{(\rr^{p+q})}^*
   \end{align*}
   and $P$\ is the stabilizer of the induced filtration of $\g$.
   The adjoint action identifies $\g_0$\ with $\co(\g_{-1},g)$.

  There is a unique (up to isomorphisms of $\g$)
  Cartan connection \[\al:\h\goesto\g\]
  which induces an isomorphism of the Euclidean spaces $(\n,\gs)$\
  and $\g/\p=\g_{-1}=(\rr^{p+q},g)$\ and satisfies
  the following two normalization conditions:
  \begin{enumerate}[(Conf.1)]
  \item \label{conformal-cond1} $\Im\kaf\subset\p$:\ its curvature has
    vanishing $\g_-$-part
  \item \label{conformal-cond2} The Ricci-type
    contraction of the $\g_0$-component of $\kaf$\ vanishes.
  \end{enumerate}

   We will now obtain maps $\Psi:K\goesto P$\ and $\al_0:\h\goesto\g$\
   such that the conformal structure $\hatgs$\ is the underlying
   structure of the homogeneous Cartan geometry corresponding to
   $(\Psi,\al_0)$.
   First note that we have a canonical embedding $\Psi$\ of $K=O(p)\times
   O(q)$\
   into $O(g)<CO(g)=G_0<P$; Being less formal, one could
   say that $\Psi=\Ad_{|K}$.
   Now $\Psi':\k\goesto\h$\ is extended to a map
   $\al_0:\h\goesto\g$\ by the obvious isometry
   of $(\n,\gs)=(\rr^p\oplus\rr^q,\frac{1}{s}\gone\oplus \frac{1}{s'}\gtwo)$\ with
   $(\g_{-1},g)=(\rr^{p+q},\gone\oplus \sgn(s')\gtwo)$:
   \begin{align*}
     {\al_0}_{|\n}:\n &\goesto\rr^{p+q}<\g,\\
     \bigl(v\oplus 0\bigr)\oplus \bigl(w\oplus 0\bigr)&\mapsto
     \bigl(\frac{1}{\sqrt{s}}v\oplus\frac{1}{\sqrt{|s'|}}w\bigr)
     \oplus 0\oplus 0.    
   \end{align*}

   It is clear that $\al_0:\h\goesto\g$\ is $K$-equivariant;
   also, $\al_0$ satisfies (C.\ref{cartcon1})\ and (C.\ref{cartcon3})\ by
   construction and thus it is a Cartan connection of type $(G,P)$.

   Since in $\o(p+1)=\rr^p\oplus\o(p)$\ the $\rr^p$-part
   brackets into $\o(p)$\ by \eqref{bracket-opplus1}\
   we also have $[\n,\n]\subset \k$; thus
   \begin{align*}
     \kaf(v_1\oplus w_1,v_2\oplus w_2)\subset\p,
   \end{align*}
   which is the first normalizing condition (Conf.\ref{conformal-cond1})\
   on the Cartan connection to this conformal structure.
   So it remains to find a map $A:\g_{-1}\goesto\g_1$\
   such that $\al=\al_0+ A\circ\al_0$\ also satisfies 
   (Conf.\ref{conformal-cond2}).

   This problem is solved by the Rho-tensor
   $A_{ij}\in L(\rr^{p+q},{(\rr^{p+q})}^*)$:
   \begin{align*}
      A_{ij}&=-\frac{1}{p+q-2}\bigl(R_{ij}-\frac{\tilde R}{2(p+q-1)}g_{ij}\bigr).
   \end{align*}
   Here $R=\kaf_0$, the $\g_0$-component of
   the curvature of $\al_0$,
   $R_{ij}={{R_{ai}}^a}_j$\ is the Ricci curvature
   and $\tilde R=g^{ij}R_{ij}$\ is the scalar curvature.
   It is a straightforward calculation that
   the Ricci curvature is
   \begin{align*}
     R_{ij}=
     s(p-1)\gone_{ij}+s'(q-1)\sgn(s')\gtwo_{ij}
   \end{align*}
   and the scalar curvature is
   \begin{align*}
     \tilde R=sp(p-1)+s'q(q-1).
   \end{align*} 

   Thus the Rho-tensor is given by
   \begin{align}\label{schouten}
     A_{ij}&=-\frac{1}{2\de}\bigl(
     (2s\De+m(s,s'))\gone_{ij}
     +(2s'\De-m(s,s'))\sgn(s')\gtwo_{ij}
     \bigr),
   \end{align}
   where
   \begin{align*}
     \de&=(p+q-1)(p+q-2);\\
     \De&=(p-1)(q-1);\\
     m(s,s')&=sp(p-1)-s'q(q-1).
   \end{align*}

   \subsection{Calculation of the curvature and the holonomy}
   Now we can calculate the curvature $\kaf$\ of the normal Cartan
   connection $\al=\al_0+A\circ\al_0$. Since, as we
   have already observed,$[\n,\n]\subset\k$, and
   since the projection of $\al$\ to $\g_0$\ vanishes,
   it is easy to see that the $\g_1$-component of $\kaf$\
   vanishes. Thus $\kaf=\kaf_{\g_0}$\ and for $X_1,X_2\in\rr^{p+q}$
   \begin{align*}
     \kaf_{\g_0}(X_1,X_2)=R(X_1,X_2)+[X_1,A(X_2)]-[X_2,A(X_1)].
   \end{align*}
   One obtains that
   ${{\kaf_{ij}}^r}_s$\ is the skew-symmetrization
   of ${{\tildka_{ij}}^r}_s$\ in the variables $i,j$, where
   \begin{align}
     {{{\tildka}_{ij}}^r}_s=
      \bigl(s-\frac{1}{\de}(m(s,s')+2s\De)\bigr)
      \deone^r_i\gone_{js}
      \notag 
       +
      \bigl(s'+\frac{1}{\de}(m(s,s')-2s'\De)\bigr)\detwo^r_i\sgn(s')\gtwo_{js}
      \notag \\ \label{kappa}
      \equghost-\frac{2\De}{\de}(s+s')
      \bigl(\deone^r_i\gtwo_{js}-\detwo^r_i\sgn(s')\gone_{js}\bigr).
   \end{align}
   \begin{thm}[Special cases]\label{thmspec}
     \begin{enumerate}
     \item      If either $p$\ or $q$\ is $1$\ or if $s'=-s$,
     $(S^p\times S^q,\gs)$\ is conformally flat.
   \item \label{spec2} If $p,q\geq 2$\ and $s'=\frac{p-1}{q-1}s$\ 
     then $(S^p\times S^q,\gs)$\ is Einstein and 
         \[\hol(S^p\times S^q,\gs)=\so(p+q+1).\]         
     \end{enumerate}
   \end{thm}
   \begin{proof}
   Of course the case $s'=-s$\ just reflects the fact that
   then $S^p\times S^q$\ is simply a twofold covering of
   the homogeneous model of conformal Cartan geometries of
   signature $(p,q)$.
   But both conformally flat cases can be seen immediately by choosing
   an orthogonal basis for $\rr\oplus\rr^{p+q}\oplus\rr$\
   and using the embeddings
   \begin{align*}
     \left(
       \begin{matrix}
       0 & -v^t \\
       v & A_1
       \end{matrix}
     \right)     
     \oplus
     \left(
       \begin{matrix}
       0 & -w^t \\
       w & A_2
       \end{matrix}
     \right) 
     \mapsto
     \left(
       \begin{matrix}
         0 & -v^t & 0 & 0 \\
         v & A_1 & 0 & 0 \\
         0 & 0 & A_2 & w \\
         0 & 0 & -w^t & 0
       \end{matrix}
     \right)
   \end{align*}
   for $s=1,s'=-1$\ resp.
   \begin{align*}
     \left(
       \begin{matrix}
       0 & -v^t \\
       v & A_1
       \end{matrix}
     \right)     
     \oplus
     \left(
       \begin{matrix}
       0 & -w \\
       w & 0
       \end{matrix}
     \right) 
     \mapsto
     \left(
       \begin{matrix}
         0 & -v^t & 0 & 0 \\
         v & A_1 & 0 & 0 \\
         0 & 0 & 0 & \frac{1}{\sqrt{|s'|}}w \\
         0 & 0 & \sgn(s')\frac{1}{\sqrt{|s'|}}w & 0
       \end{matrix}
     \right)
   \end{align*}
   for $q=1,s=1$ and $s'$\ of arbitrary signature.

\noindent
   To see the second claim 
       we first notice that $\hol(\al)$\ contains $\so(\rr^n,g)$:
       Since \[\frac{2\De}{\de}(s+s')\not=0\] the image of $\kaf$\ contains all matrices
       of the form
       \begin{align*}
         \left(
         \begin{matrix}
           0 & B \\
           -\sgn(s') B^t & 0
         \end{matrix}
         \right)\in\so(\rr^n,g).
       \end{align*}
       But since
       \begin{align*}
       \begin{matrix}
         \al\left(
         \left(
         \begin{matrix}
           0 & 0 \\
           0 & A_1 
         \end{matrix}
         \right)
         \oplus
         \left(
         \begin{matrix}
           0 & 0 \\
           0 & A_2 
         \end{matrix}
         \right)
         \right)
         =
         \left(
         \begin{matrix}
           A_1 & 0 \\
           0  & A_2
         \end{matrix}
         \right)\in\so(\rr^n,g)
       \end{matrix}
       \end{align*}
       and $\hol(\al)$\ is the $(\h,\al)$-module generated by the image of $\kaf$, 
       simplicity of $\so(\rr^n,g)$\ shows that indeed $\so(\rr^n,g)\subset\hol(\al)$.

\noindent
       Now the condition
       \begin{align*}
         s'=\frac{p-1}{q-1}s
       \end{align*}
       means exactly that the Ricci curvature $R_{ij}$\ is a multiple
       of $g$, and thus $(S^p\times S^q,g_{(s,s')})$\ is Einstein.
       Then by \eqref{schouten}\  also $A$\ is a multiple of $g$, and,
       explicitly: $A=rg$\ with
       \begin{align*}
         r=-\frac{p-1}{2(p+q-1)}.
       \end{align*}
       Thus
       \begin{align*}
         \al( (v\oplus A_1)\oplus (w\oplus A_2)=(v\oplus w)\oplus(A_1\oplus A_2)\oplus
         r(v\oplus w)^t,
       \end{align*}
       and we see that the smallest Lie subalgebra of $\g$\ containing $\Im\kaf$
       and being invariant under $(\h,\al)$\ consists of matrices of the form
       \begin{align*}
         \left(
         \begin{matrix}
           0 & r X^t & 0 \\
           X & A & -r X \\
           0 & -X^t & 0
         \end{matrix}
         \right),
       \end{align*}
       which proves the claim since $r<0$.       
   \end{proof}
   When both $p$\ and $q$\ are at least two a generic ratio of radii
   $(s,s')$\ of arbitrary signature yields full holonomy, which is straightforward to check. A different argument for this case can be found in \cite{armstrong-conformal}.
   
   \begin{thm}[The generic case]
     If $p,q\geq 2$, then for $s'\not\in\{-s,\frac{p-1}{q-1}s\}$,
     \begin{align*}
     \hol(S^p\times S^q,\gs)&=\g=\so(p+q+1,1) \ \mathrm{for}\ s'>0\ \mathrm{resp.}\\
     \hol(S^p\times S^q,\gs)&=\g=\so(p+1,q+1)\ \mathrm{for}\ s'<0.
     \end{align*}
   \end{thm}

   \begin{rema}
       The treatment of the holonomy of homogeneous parabolic geometries other than conformal
       structures is closely parallel, the only additional problem
       which appears is that there are no longer general
formulas for the prolongation of the given geometric data to the
corresponding Cartan geometries. To see how this problem boils down to basic representation theory see e.g. \cite{cap-2006-cr}\ or \cite{mrh-dip}\ for explicit
examples of prolongations in the realm of CR-structures.
   \end{rema}

\end{document}